\documentclass[11pt]{article}
\usepackage{amssymb}
\usepackage{latexsym}
\title{On identities in Hom-Malcev algebras}
\author{{\it A. Nourou Issa} \\ {\it D\'epartement de Math\'ematiques, Universit\'e d'Abomey-Calavi} \\ {\it 01 BP 4521 Cotonou 01, Benin}. \\ {\it E-mail: woraniss@yahoo.fr}}
\date{}
\begin{document}
\maketitle
\begin{abstract}
In a Hom-algebra an identity, equivalent to the Hom-Malcev identity, is found.
\par
2000 MSC: 17A20, 17A30.
\end{abstract}
\section{Introduction and statement of results}
Hom-Lie algebras were introduced in [3] as a tool in understanding the structure of some q-deformations of the Witt and the Virasoro algebras. Since then, the theory of Hom-type algebras began an intensive development (see, e.g., [2], [4], [6], [7], [8], [12], [13], [14], [15]). Hom-type algebras are defined by twisting the defining identities of some well-known algebras by a linear self-map, and when this twisting map is the identity map, one recovers the original type of considered algebras.
\par
In this setting, a Hom-type generalization of Malcev algebras (called Hom-Malcev algebras) is defined by D. Yau in [15]. Recall that a {\it Malcev algebra} is a nonassociative algebra $(A, \cdot )$, where the binary operation ``$\cdot$'' is anti-commutative, such that the identity \\
\par
$J(x,y,x \cdot z) = J(x,y,z) \cdot x$ \hfill (1.1) \\
\\
holds for all $x,y,z \in A$ (here $J(x,y,z)$ denotes the Jacobian, i.e. $J(x,y,z) = xy \cdot z + yz \cdot x + zx \cdot y$). The identity (1.1) is known as the {\it Malcev identity}. Malcev algebras were introduced by A.I. Mal'tsev [9] (calling them Moufang-Lie algebras) as tangent algebras to local smooth loops, generalizing in this way a result in Lie theory stating that a Lie algebra is a tangent algebra to a local Lie group (in fact, Lie algebras are special case of Malcev algebras). Another approach to Malcev algebras is the one from alternative algebras: every alternative algebra is Malcev-admissible [9]. So one could say that the algebraic theory of Malcev algebras started from Malcev-admissibility of algebras. The foundations of the algebraic theory of Malcev algebras go back to E. Kleinfeld [5], A.A. Sagle [10] and, as mentioned in [10], to A.A. Albert and L.J. Paige. Some twisting of the Malcev identity (1.1) along any algebra self-map $\alpha$ of $A$ gives rise to the notion of a Hom-Malcev algebra $(A, \cdot , \alpha)$ ([15]; see definitions in section 2). Properties and constructions of Hom-Malcev algebras, as well as the relationships between these Hom-algebras and Hom-alternative or Hom-Jordan algebras are investigated in [15]. In particular, it is shown that a Malcev algebra can be twisted into a Hom-Malcev algebra and that Hom-alternative algebras are Hom-Malcev admissible.
\par
In [15], as for Malcev algebras (see [10], [11]), equivalent defining identities of a Hom-Malcev algebra are given. In this note, we mention another identity in a Hom-Malcev algebra that is equivalent to the ones found in [15]. Specifically, we shall prove the following \\
\par
{\bf Theorem.} {\it Let $(A, \cdot , \alpha)$ be a Hom-Malcev algebra. Then the identity} \\
\par
$J_{\alpha}(w \cdot x, \alpha (y), \alpha (z)) = J_{\alpha}(w,y,z) \cdot {\alpha}^{2}(x) + {\alpha}^{2}(w) \cdot J_{\alpha}(x,y,z)$ 
\par
\hspace{3.5cm} $- 2J_{\alpha}(y \cdot z, \alpha (w), \alpha (x))$  \hfill (1.2) \\
\\
{\it holds for all $w,x,y,z$ in $A$, where $J_{\alpha}(x,y,z) = xy \cdot \alpha (z) + yz \cdot \alpha (x) + zx \cdot \alpha (y)$. Moreover, in any anti-commutative Hom-algebra $(A, \cdot , \alpha)$, the identity (1.2) is equivalent to the Hom-Malcev identity} \\
\par
$J_{\alpha}(\alpha (x), \alpha (y), x \cdot z) = J_{\alpha}(x,y,z) \cdot {\alpha}^{2}(x)$ \hfill (1.3) \\
\\
{\it for all $x,y,z$ in $A$}. \\
\par
Observe that when $\alpha = Id$ (the identity map) in (1.3), then (1.3) is (1.1) i.e. the Hom-Malcev algebra $(A, \cdot , \alpha)$ reduces to the Malcev algebra $(A, \cdot )$ (see [15]).
\par
In section 2 some instrumental lemmas are proved. Some results in these lemmas are a kind of the Hom-version of similar results by E. Kleinfeld [5] in case of Malcev algebras. The section 3 is devoted to the proof of the theorem.
\par
Throughout this note we work over a ground field $\mathbb K$ of characteristic 0.
\section{Definitions. Preliminary results}
In this section we recall useful notions on Hom-algebras ([8], [12], [13], [15]), as well as the one of a Hom-Malcev algebra [15]. In [5], using an analogue of the Bruck-Kleinfeld function, an identity (see identity (6) in [5]) characterizing Malcev algebras is found. This identity is used in [10] to derive further identities for Malcev algebras (see [10], Proposition 2.23). The main result of this section (Lemma 2.7) proves that the Hom-version of the identity (6) of [5] holds in any Hom-Malcev algebra. \\
\par
{\bf Definition 2.1.} A {\it multiplicative Hom-algebra} is a triple $(A, \mu, \alpha)$ , in which $A$ is a $\mathbb K$-module, $\mu : A \times A \rightarrow A$ is a bilinear map (the binary operation), and  $\alpha : A \rightarrow A$ is a linear map (the twisting map) such that $\alpha$ is an endomorphism of $(A, \mu)$. The Hom-algebra $(A, \mu, \alpha)$ is said {\it anticommutative} if the operation $\mu$ is skew-symmetric, i.e. $\mu (x,y) = - \mu (y,x)$, for all $x,y \in A$. \\
\par
In the rest of this paper, we will use the abbreviation $x \cdot y = \mu (x,y)$ in a Hom-algebra $(A, \mu, \alpha)$. \\
\par
{\bf Remark.} The multiplicativity of the twisting map is not necessary in the definition of a Hom-algebra (see, e.g., [6], [8]). The multiplicativity is included here for convenience. \\
\par
{\bf Definition 2.2.} Let $(A, \cdot, \alpha)$ be an anticommutative Hom-algebra.
\par
(i) The {\it Hom-Jacobian} ([8]) of $(A, \cdot, \alpha)$ is the trilinear map $ J_{\alpha}(x,y,z)$ on A defined by $J_{\alpha}(x,y,z) = xy \cdot \alpha (z) + yz \cdot \alpha (x) + zx \cdot \alpha (y)$.
\par
(ii) $(A, \cdot, \alpha)$ is called a {\it Hom-Lie algebra} ([3]) if the {\it Hom-Jacobi identity} $ J_{\alpha}(x,y,z) = 0$ holds in $(A, \cdot, \alpha)$. \\
\par
{\bf Definition 2.3.} ([15]) A {\it Hom-Malcev algebra} is an anticommutative algebra $(A, \cdot, \alpha)$ such that the {\it Hom-Malcev identity} (see (1.3)) \\
\par
$J_{\alpha}(\alpha (x), \alpha (y), x \cdot z) = J_{\alpha}(x,y,z) \cdot {\alpha}^{2}(x)$ \\
\\
holds in $(A, \cdot, \alpha)$. \\
\par
{\bf Remark.} When $\alpha = Id$, then the Hom-Jacobi identity reduces to the usual Jacobi identity $J(x,y,z) := xy \cdot z + yz \cdot x + zx \cdot y = 0$, i.e. the Hom-Lie algebra $(A, \cdot, \alpha)$ reduces to the Lie algebra $(A, \cdot )$. Likewise, when $\alpha = Id$, the Hom-Malcev identity reduces to the Malcev identity (1.1), i.e. the Hom-Malcev algebra $(A, \cdot, \alpha)$ reduces to the Malcev algebra $(A, \cdot )$. \\
\par
The following simple lemma holds in any anticommutative Hom-algebra. \\
\par
{\bf Lemma 2.4.} {\it In any anticommutative Hom-algebra $(A, \cdot, \alpha)$ the following holds:
\par
(i) $J_{\alpha}(x,y,z)$ is skew-symmetric in its three variables.
\par
(ii) ${\alpha}^{2}(w) \cdot J_{\alpha}(x,y,z) - {\alpha}^{2}(x) \cdot J_{\alpha}(y,z,w) + {\alpha}^{2}(y) \cdot J_{\alpha}(z,w,x)$ 
\par 
\hspace{1cm} $ - {\alpha}^{2}(z) \cdot J_{\alpha}(w,x,y)$
\par
\hspace{1cm} $= J_{\alpha}(w \cdot x, \alpha (y), \alpha (z)) + J_{\alpha}(y \cdot z, \alpha (w), \alpha (x)) + J_{\alpha}(w \cdot y, \alpha (z), \alpha (x))$
\par
\hspace{1cm} $ + J_{\alpha}(z \cdot x, \alpha (w), \alpha (y)) - J_{\alpha}(z \cdot w, \alpha (x), \alpha (y)) - J_{\alpha}(x \cdot y, \alpha (z), \alpha (w))$, \\
for all $w,x,y, z$ in $A$}. \\
\par
{\bf Proof.} The skew-symmetry of $J_{\alpha}(x,y,z)$ in $w,x,y, z$ follows from the skew-symmetry of the operation ``$\cdot$''.
\par
Expanding the expression in the left-hand side of (ii) and then rearranging terms, we get (by the skew-symmetry of ``$\cdot$'')
\begin{eqnarray}
& &{\alpha}^{2}(w) \cdot J_{\alpha}(x,y,z) - {\alpha}^{2}(x) \cdot J_{\alpha}(y,z,w) + {\alpha}^{2}(y) \cdot J_{\alpha}(z,w,x) \nonumber \\ 
& & - {\alpha}^{2}(z) \cdot J_{\alpha}(w,x,y) \nonumber \\
& & = - {\alpha}^{2}(z) \cdot (wx \cdot \alpha (y)) + {\alpha}^{2}(y) \cdot (wx \cdot \alpha (z))\nonumber \\
& & - {\alpha}^{2}(x) \cdot (yz \cdot \alpha (w)) + {\alpha}^{2}(w) \cdot (yz \cdot \alpha (x))\nonumber \\
& & - {\alpha}^{2}(x) \cdot (wy \cdot \alpha (z)) - {\alpha}^{2}(z) \cdot (yw \cdot \alpha (x))\nonumber \\
& & + {\alpha}^{2}(w) \cdot (zx \cdot \alpha (y)) + {\alpha}^{2}(y) \cdot (xz \cdot \alpha (w))\nonumber \\
& & - {\alpha}^{2}(x) \cdot (zw \cdot \alpha (y)) + {\alpha}^{2}(y) \cdot (zw \cdot \alpha (x))\nonumber \\
& & + {\alpha}^{2}(w) \cdot (xy \cdot \alpha (z)) - {\alpha}^{2}(z) \cdot (xy \cdot \alpha (w)).\nonumber 
\end{eqnarray}
Next, adding and subtracting $\alpha (yz) \cdot \alpha (wx)$ (resp. $\alpha (wx) \cdot \alpha (yz)$, \\ $\alpha (zx) \cdot \alpha (wy)$, $\alpha (wy) \cdot \alpha (zx)$, $\alpha (xy) \cdot \alpha (zw)$ and $\alpha (zw) \cdot \alpha (xy)$) in the first (resp. second, third, fourth, fifth, and sixth) line of the right-hand side expression in the last equality above, we come to the equality (ii) of the lemma. \hfill $\square$ \\
\par
In a Hom-Malcev $(A, \cdot, \alpha)$ we define the multilinear map $G$ by \\
\par
$ G(w,x,y,z) =  J_{\alpha}(w \cdot x, \alpha (y), \alpha (z)) - {\alpha}^{2}(x) \cdot J_{\alpha}(w,y,z)$ 
\par
\hspace{2cm} $ - J_{\alpha}(x,y,z) \cdot {\alpha}^{2}(w)$ \hfill (2.1) \\
\\
for all $w,x,y,z$ in $A$. \\
\par
{\bf Remark.} (1) If $\alpha = Id$ in (2.1), then $ G(w,x,y,z)$ reduces to the function $ f(w,x,y,z)$ defined in [5] which in turn is a variation of the Bruck-Kleinfeld function defined in [1].
\par
(2) If in (2.1) replace $J_{\alpha}(t,u,v)$ with the {\it Hom-associator} [8] $as(t,u,v)$, then one recovers the Hom-Bruck-Kleinfeld function defined in [15]. \\
\par
{\bf Lemma 2.5.} {\it In a Hom-Malcev algebra $(A, \cdot, \alpha)$ the function $G(w,x,y,z)$ defined by (2.1) is skew-symmetric in its four variables}. \\
\par
{\bf Proof.} From the skew-symmetry of ``$\cdot$'' and $J_{\alpha}(t,u,v)$ (see Lemma 2.4(i)) it clearly follows that
\par
$G(x,w,y,z) = - G(w,x,y,z)$
\par
$G(w,x,z,y) = - G(w,x,y,z)$. \\
Next, using the skew-symmetry of $J_{\alpha}(t,u,v)$, 
\begin{eqnarray}
G(y,x,y,z) &=& J_{\alpha}(y \cdot x, \alpha (y), \alpha (z)) - J_{\alpha}(x,y,z) \cdot {\alpha}^{2}(y)\nonumber \\
&=& J_{\alpha}(\alpha (y), \alpha (z), y \cdot x) - J_{\alpha}(y,z,x) \cdot {\alpha}^{2}(y)\nonumber \\
&=& J_{\alpha}(y,z,x) \cdot {\alpha}^{2}(y) - J_{\alpha}(y,z,x) \cdot {\alpha}^{2}(y) \;\; \mbox{(by (1.3))}
\nonumber \\
&=& 0. \nonumber
\end{eqnarray}
Likewise, one checks that $G(w,y,y,z) = 0$. This suffices to prove the skew-symmetry of $G(w,x,y,z)$ in its variables. \hfill $\square$ \\
\par
As we shall see below, the following lemma is a consequence of the definition of $G(w,x,y,z)$ and the skew-symmetry of $J_{\alpha}(t,u,v)$ and $G(w,x,y,z)$. \\
\par
{\bf Lemma 2.6.} {\it Let $(A, \cdot, \alpha)$ be a Hom-Malcev . Then} \\
\par
$ J_{\alpha}(w \cdot x, \alpha (y), \alpha (z)) + J_{\alpha}(x \cdot y, \alpha (z), \alpha (w)) + J_{\alpha}(y \cdot z, \alpha (w), \alpha (x))$ 
\par
$+ J_{\alpha}(z \cdot w, \alpha (x), \alpha (y)) = 0$; \hfill (2.2) \\
\par
$ 2 G(w,y,y,z) - {\alpha}^{2}(w) \cdot J_{\alpha}(x,y,z) + {\alpha}^{2}(x) \cdot J_{\alpha}(w,y,z) - {\alpha}^{2}(y) \cdot J_{\alpha}(z,w,x) $ 
\par
$+ {\alpha}^{2}(z) \cdot J_{\alpha}(w,x,y) = J_{\alpha}(w \cdot x, \alpha (y), \alpha (z)) + J_{\alpha}(y \cdot z, \alpha (w), \alpha (x))$, \hfill (2.3)\\
\\
{\it for all $w,x,y,z$ in A}. \\
\par
{\bf Proof.} From the definition of $G(w,x,y,z)$ (see (2.1)) we have 
\par
$ J_{\alpha}(w \cdot x, \alpha (y), \alpha (z)) = G(w,x,y,z) + {\alpha}^{2}(x) \cdot J_{\alpha}(w,y,z) + J_{\alpha}(x,y,z) \cdot {\alpha}^{2}(w)$,
\par
$ J_{\alpha}(x \cdot y, \alpha (z), \alpha (w)) = G(x,y,z,w) +  {\alpha}^{2}(y) \cdot J_{\alpha}(x,z,w) + J_{\alpha}(y,z,w) \cdot {\alpha}^{2}(x)$,
\par
$ J_{\alpha}(y \cdot z, \alpha (w), \alpha (x)) = G(y,z,w,x) +  {\alpha}^{2}(z) \cdot J_{\alpha}(y,w,x) + J_{\alpha}(z,w,x) \cdot {\alpha}^{2}(y)$,
\par
$ J_{\alpha}(z \cdot w, \alpha (x), \alpha (y)) = G(z,w,x,y) + {\alpha}^{2}(w) \cdot J_{\alpha}(z,x,y) + J_{\alpha}(w,x,y) \cdot {\alpha}^{2}(z)$. \\
Therefore, by the skew-symmetry of ``$\cdot$'', $J_{\alpha}(x,y,z)$ and $G(w,x,y,z)$, we get
\par
$ J_{\alpha}(w \cdot x, \alpha (y), \alpha (z)) + J_{\alpha}(x \cdot y, \alpha (z), \alpha (w)) + J_{\alpha}(y \cdot z, \alpha (w), \alpha (x))$ 
\par
$ + J_{\alpha}(z \cdot w, \alpha (x), \alpha (y))$
\par
$= G(w,x,y,z) + G(x,y,z,w) + G(y,z,w,x) + G(z,w,x,y)$ 
\par
$= G(w,x,y,z) - G(w,x,y,z) + G(y,z,w,x) - G(y,z,w,x)$
\par
$=0$, \\
which proves (2.2).
\par
Next, again from the expression of $G(w,x,y,z)$,
\par
$ J_{\alpha}(w \cdot x, \alpha (y), \alpha (z)) + J_{\alpha}(y \cdot z, \alpha (w), \alpha (x))$
\par
$= [G(w,x,y,z) + {\alpha}^{2}(x) \cdot J_{\alpha}(w,y,z) + J_{\alpha}(x,y,z) \cdot {\alpha}^{2}(w)] $
\par
$+ [G(y,z,w,x) + {\alpha}^{2}(z) \cdot J_{\alpha}(y,w,x) + J_{\alpha}(z,w,x) \cdot {\alpha}^{2}(y)]$
\par
$= 2 G(w,x,y,z) - {\alpha}^{2}(w) \cdot J_{\alpha}(x,y,z) + {\alpha}^{2}(x) \cdot J_{\alpha}(y,z,w) - {\alpha}^{2}(y) \cdot J_{\alpha}(z,w,x)$
\par
$+ {\alpha}^{2}(z) \cdot J_{\alpha}(w,x,y)$ \\
so that we get (2.3). \hfill $\square$ \\
\par
From Lemma 2.5 and Lemma 2.6, we get the following expression of $G(w,x,y,z)$. \\
\par
{\bf Lemma 2.7.} {\it Let $(A, \cdot, \alpha)$ be a Hom-Malcev . Then} \\
\par
$G(w,x,y,z) = 2[J_{\alpha}(w \cdot x, \alpha (y), \alpha (z)) + J_{\alpha}(y \cdot z, \alpha (w), \alpha (x))] $ \hfill (2.4) \\
\\
{\it for all $w,x,y,z$ in $A$}. \\
\par
{\bf Proof.} Set $ g(w,x,y,z) = J_{\alpha}(w \cdot x, \alpha (y), \alpha (z)) + J_{\alpha}(x \cdot y, \alpha (z), \alpha (w)) + J_{\alpha}(y \cdot z, \alpha (w), \alpha (x)) + J_{\alpha}(z \cdot w, \alpha (x), \alpha (y))$. Then (2.2) says that $ g(w,x,y,z) = 0$ for all $ w,x,y,z $ in $A$. Now, by adding $ g(w,x,y,z) - g(x,w,y,z)$ to the right-hand side of Lemma 2.4(ii), we get
\begin{eqnarray}
{\alpha}^{2}(w) \cdot J_{\alpha}(x,y,z) &-& {\alpha}^{2}(x) \cdot J_{\alpha}(y,z,w) \nonumber \\
&+& {\alpha}^{2}(y) \cdot J_{\alpha}(z,w,x) - {\alpha}^{2}(z) \cdot J_{\alpha}(w,x,y) \nonumber \\
&=& J_{\alpha}(w \cdot x, \alpha (y), \alpha (z)) + J_{\alpha}(y \cdot z, \alpha (w), \alpha (x))\nonumber \\
&+& J_{\alpha}(w \cdot y, \alpha (z), \alpha (x)) + J_{\alpha}(z \cdot x, \alpha (w), \alpha (y))\nonumber \\
&-& J_{\alpha}(z \cdot w, \alpha (x), \alpha (y)) - J_{\alpha}(x \cdot y, \alpha (z), \alpha (w))\nonumber \\
&+& J_{\alpha}(w \cdot x, \alpha (y), \alpha (z)) + J_{\alpha}(x \cdot y, \alpha (z), \alpha (w))\nonumber \\
&+& J_{\alpha}(y \cdot z, \alpha (w), \alpha (x)) + J_{\alpha}(z \cdot w, \alpha (x), \alpha (y))\nonumber \\
&-& J_{\alpha}(x \cdot w, \alpha (y), \alpha (z)) - J_{\alpha}(w \cdot y, \alpha (z), \alpha (x))\nonumber \\
&-& J_{\alpha}(y \cdot z, \alpha (x), \alpha (w)) - J_{\alpha}(z \cdot x, \alpha (w), \alpha (y))\nonumber \\
&=& 3J_{\alpha}(w \cdot x, \alpha (y), \alpha (z)) + 3J_{\alpha}(y \cdot z, \alpha (w), \alpha (x)) \nonumber
\end{eqnarray}
i.e.
\begin{eqnarray}
{\alpha}^{2}(w) \cdot J_{\alpha}(x,y,z) &-& {\alpha}^{2}(x) \cdot J_{\alpha}(y,z,w) \nonumber \\
&+& {\alpha}^{2}(y) \cdot J_{\alpha}(z,w,x) - {\alpha}^{2}(z) \cdot J_{\alpha}(w,x,y) \nonumber \\
&=& 3[J_{\alpha}(w \cdot x, \alpha (y), \alpha (z)) + J_{\alpha}(y \cdot z, \alpha (w), \alpha (x))]. \;\; \mbox{(2.5)} \nonumber
\end{eqnarray}
Next, adding (2.3) and (2.5) together, we get
\begin{eqnarray}
2G(w,x,y,z) &-& {\alpha}^{2}(w) \cdot J_{\alpha}(x,y,z)+{\alpha}^{2}(x) \cdot J_{\alpha}(y,z,w) \nonumber \\
&-& {\alpha}^{2}(y) \cdot J_{\alpha}(z,w,x) + {\alpha}^{2}(z) \cdot J_{\alpha}(w,x,y) \nonumber \\
&+& {\alpha}^{2}(w) \cdot J_{\alpha}(x,y,z) - {\alpha}^{2}(x) \cdot J_{\alpha}(y,z,w) \nonumber \\
&+& {\alpha}^{2}(y) \cdot J_{\alpha}(z,w,x) - {\alpha}^{2}(z) \cdot J_{\alpha}(w,x,y) \nonumber \\
&=& J_{\alpha}(w \cdot x, \alpha (y), \alpha (z)) + J_{\alpha}(y \cdot z, \alpha (w), \alpha (x))\nonumber \\
&+& 3[J_{\alpha}(w \cdot x, \alpha (y), \alpha (z))+J_{\alpha}(y \cdot z, \alpha (w), \alpha (x))] \nonumber
\end{eqnarray}
i.e. 
\par
$ 2G(w,x,y,z) = 4[J_{\alpha}(w \cdot x, \alpha (y), \alpha (z))+J_{\alpha}(y \cdot z, \alpha (w), \alpha (x))]$ \\
and (2.4) follows. \hfill $\square$
\section{Proof}
Relaying on the lemmas of section 2, we are now in position to prove the theorem. \\
\par
{\bf Proof of the theorem.} First we establish the identity (1.2) in a Hom-Malcev algebra. We may write (2.1) in an equivalent form:\\
\par
$J_{\alpha}(w \cdot x, \alpha (y), \alpha (z)) =  {\alpha}^{2}(x) \cdot J_{\alpha}(w,y,z)
+ J_{\alpha}(x,y,z) \cdot {\alpha}^{2}(w)$
\par
\hspace{3.5cm} $+ G(w,x,y,z)$. \hfill (2.6) \\
\\
Now in (2.6), replace $G(w,x,y,z)$ with its expression from (2.4) to get
$-J_{\alpha}(w \cdot x, \alpha (y), \alpha (z)) =  {\alpha}^{2}(x) \cdot J_{\alpha}(w,y,z)
+ J_{\alpha}(x,y,z) \cdot {\alpha}^{2}(w)$
\par
\hspace{3.5cm} $ +2 J_{\alpha}(y \cdot z, \alpha (w), \alpha (x))$, \\
which leads to (1.2).
\par
Now, we proceed to prove the equivalence of (1.2) with (1.3) in an anticommutative Hom-Malcev algebra. First assume (1.3). Then Lemmas 2.4, 2.5, 2.6, and 2.7 imply that (1.2) holds in any Hom-Malcev algebra.
\par
Conversely, assume (1.2). Then, setting $w = y$ in (1.2), we get, by the skew-symmetry of $J_{\alpha}(x,y,z)$, \\
\par
$ J_{\alpha}(y \cdot x, \alpha (y), \alpha (z)) = {\alpha}^{2}(y) \cdot J_{\alpha}(y,z,x) - 2 J_{\alpha}(\alpha (y), \alpha (x), y \cdot z)$. \hfill (2.7) \\
\\
Now, the permutation of $z$ with $x$ in (2.7) gives \\
$ J_{\alpha}(y \cdot z, \alpha (y), \alpha (x)) = {\alpha}^{2}(y) \cdot J_{\alpha}(y,x,z) - 2 J_{\alpha}(\alpha (y), \alpha (z), y \cdot x)$, \\
i.e. \\
$ 2J_{\alpha}(y \cdot z, \alpha (y), \alpha (x)) = -2{\alpha}^{2}(y) \cdot J_{\alpha}(y,z,x) - 4 J_{\alpha}(\alpha (y), \alpha (z), y \cdot x)$, \\
or \\
\par
$  4 J_{\alpha}(\alpha (y), \alpha (z), y \cdot x)= -2{\alpha}^{2}(y) \cdot J_{\alpha}(y,z,x) - 2J_{\alpha}(y \cdot z, \alpha (y), \alpha (x))$.  \hfill (2.8) \\
\\
Next, the subtraction of (2.8) from (2.7) gives (keeping in mind the skew-symmetry of $J_{\alpha}(x,y,z)$) \\
$-3J_{\alpha}(\alpha (y), \alpha (z), y \cdot x)= 3{\alpha}^{2}(y) \cdot J_{\alpha}(y,z,x)$ \\
i.e. \\
$J_{\alpha}(\alpha (y), \alpha (z), y \cdot x)=J_{\alpha}(y,z,x) \cdot {\alpha}^{2}(y)$ \\
so that we get (1.3). \hfill $\square$ \\
\par
{\bf Remark.} If set $\alpha = Id$, then the identity (1.2) (resp. (1.3)) reduces to the identity (2.26) (resp. (2.4)) of [10]. The equivalence of (2.4) and (2.26) of [10] could be deduced from the works [10] and [11].

\end{document}